\documentclass[12pt]{amsart}
\usepackage[cp850]{inputenc}
\usepackage{amssymb}
\usepackage{amsmath}
\usepackage{latexsym}

\newcommand\mylabel[1]{\label{#1}}			
\newcommand\comment[1]{}				

\def\bbR{\mathbb{R}}

\def\s{*}
\def\nm{\!\nmid\!}
\newcommand\cA{{\mathcal A}}
\newcommand\cP{{\mathcal P}}

\def\qedsymbol{\blacksquare}
\newcommand\eset{\varnothing}

\begin{document}
\setlength{\parindent}{0pt}
\setlength{\parskip}{0.4cm}

\newtheorem{theorem}{Theorem}[section] 
\newtheorem{lemma}[theorem]{Lemma} 
\newtheorem{corollary}[theorem]{Corollary} 
\newtheorem{proposition}[theorem]{Proposition} 

\theoremstyle{remark}
\newtheorem{example}{Example}[section]

\begin{center}

\Large{\bf A unifying generalization of Sperner's theorem}
\normalsize

{\sc Matthias Beck, Xueqin Wang, and Thomas Zaslavsky\footnote{Research supported by National Science Foundation grant DMS-70729.}}

{\sc State University of New York at Binghamton}

{\tt matthias@math.binghamton.edu \\
    xwang@math.binghamton.edu \\ 
    zaslav@math.binghamton.edu }


\end{center}
\bigskip\bigskip

{\it Abstract:} 
Sperner's bound on the size of an antichain in the lattice $\cP(S)$ of subsets of a finite set $S$ has been generalized in three different directions: 
by Erd\H{o}s to subsets of $\cP(S)$ in which chains contain at most $r$ elements; by Meshalkin to certain classes of compositions of $S$; 
by Griggs, Stahl, and Trotter through replacing the antichains by certain sets of pairs of disjoint elements of $\cP(S)$. 
We unify Erd\H{o}s's, Meshalkin's, and Griggs--Stahl--Trotter's inequalities with a common generalization.  
We similarly unify their accompanying LYM inequalities. Our bounds do not in general appear to be the best possible. 
\normalsize



{\it Running head}:  A unifying Sperner generalization

{\it Address for editorial correspondence}:\\
    Matthias Beck\\
    Department of Mathematical Sciences\\
    State University of New York\\
    Binghamton, NY 13902-6000 \\
    U.S.A.

\vfill
\pagebreak

\section{Sperner-type theorems} \mylabel{intro}
Let $ S $ be a finite set with $ n $ elements. In the lattice $ \cP(S) $ of all subsets of $ S $ one tries to estimate the size of a subset with certain characteristics.  The most famous 
such estimate concerns {\bf antichains}, that is, subsets of $ \cP(S) $ 
in which any two elements are incomparable.  We let $\lfloor x \rfloor$ denote the greatest integer $\leq x$ and $\lceil x \rceil$ the least integer $\geq x$.
\begin{theorem}[Sperner \cite{sperner}]\mylabel{sperner} 
Suppose $ A_1, \dots , A_m \subseteq S $ such that $ A_k \not\subseteq A_j $ for $ k \neq j $. Then 
$ m \leq \binom{n}{ \lfloor n/2 \rfloor } $.
Furthermore, this bound can be attained for any $n$. \end{theorem} 

Sperner's theorem has been generalized in many different directions.  Here are three: 
Erd\H{o}s extended Sperner's inequality to subsets of $\cP(S)$ in which chains contain at most $r$ elements. 
Meshalkin proved a Sperner-like inequality for families of compositions of $S$ into a fixed number of parts, in which the sets in each part constitute an antichain. 
Finally, Griggs, Stahl, and Trotter extended Sperner's theorem by replacing the antichains by sets of pairs of disjoint elements of $\cP(S)$ satisfying an intersection condition. 
In this paper we unify Erd\H{o}s's, Meshalkin's, and the Griggs--Stahl--Trotter inequalities in a single generalization. 
However, except in special cases (among which are generalizations of the known bounds), our bounds 
are not the best possible. 

For a precise statement of Erd\H{o}s's generalization, call a subset of $\cP(S)$ 
{\bf $r$-chain-free} if its chains (i.e., linearly ordered subsets) contain no more than $r$ 
elements; that is, no chain has length $r$.\footnote{The term ``$r$-family'' or ``$k$-family'', depending on the name 
of the forbidden length, has been used in the past, but we think it is time for a distinctive name.} 
In particular, an antichain is $1$-chain-free. The generalization of Theorem \ref{sperner} to $r$-chain-free families is 

\begin{theorem}[Erd\H{o}s \cite{erdos}]\mylabel{erdos} 
Suppose $ \left\{ A_1, \dots , A_m \right\} \subseteq \cP(S) $ contains no chains with $r+1$ elements. Then $m$ 
is bounded by the sum of the $r$ largest binomial coefficients $\binom{n}{k}$, $0 \leq k \leq n$.  The bound  is attainable for every $n$ and $r$. \end{theorem} 

Note that for $ r=1 $, we obtain Sperner's theorem.

Going in a different direction, Sperner's inequality can be generalized to certain ordered weak partitions of $ S $. 
We define a {\bf weak composition of $S$ into $p$ parts} as an ordered $p$-tuple $ \left( A_1, \dots , A_p \right) $ of sets $A_k$, possibly void, such that 
$ A_1, \dots , A_p $ are pairwise disjoint and $ A_1 \cup \dots \cup A_p = S $. 
A Sperner-like inequality suitable for this setting was proposed by Sevast'yanov and proved by Meshalkin (see \cite{meshalkin}).  By a {\bf $p$-multinomial coefficient for $n$} we mean a multinomial coefficient $\binom{n}{a_1, \dots, a_p}$, where $a_i \geq 0$ and $a_1 + \cdots + a_p = n$.  Let $[p] := \{1,2,\ldots,p\}$.

\begin{theorem}[Meshalkin]\mylabel{meshalkin} 
Let $p \geq 2$.  Suppose $ \left( A_{j1} , \dots , A_{jp} \right) $ for $j=1,\ldots,m$ are different weak compositions of $S$ into $p$ parts 
such that for all $ k \in [p] $ the set $ \left\{ A_{jk} : 1 \leq j \leq m \right\} $ (ignoring repetition) forms an antichain.
 Then $m$ is bounded by the largest $p$-multinomial coefficient for $n$.  Furthermore, the bound is attainable for every $n$ and $p$.
\end{theorem} 

This largest multinomial coefficient can be written explicitly as
\[ 
\frac{ n! }{\big( ( \big\lfloor \frac{n}{p} \big\rfloor + 1 )! \big)^\rho \big( \big\lfloor \frac{n}{p} \big\rfloor !\big)^{ p-\rho } } \ , 
\] 
where $ \rho = n - p \big\lfloor \frac{n}{p} \big\rfloor $.  
To see why Meshalkin's inequality generalizes Sperner's Theorem, suppose $ A_1, \dots , A_m \subseteq S $ form an 
antichain. Then $ S-A_1, \dots , S-A_m $ also form an antichain. Hence the $m$ weak compositions $ \left( A_j , S-A_j \right) $ of $S$ 
into two parts satisfy Meshalkin's conditions and Sperner's inequality follows. 

Yet another generalization of Sperner's Theorem is 

{\bf Theorem \ref{gst}$'$ (Griggs--Stahl--Trotter} \cite{stahl}{\bf )} \quad
{\it Suppose $ \left\{A_{j0} , \dots , A_{jq} \right\} $ are $m$ different chains in $\cP(S)$ such that $ A_{ji} \not\subseteq A_{kl} $ for all $i$ and $l$ and all $j \neq k$.  Then $ m \leq \binom{ n-q }{ \lfloor (n-q)/2 \rfloor } $.  Furthermore, this bound can be attained for all $n$ and $q$. }
\medskip

An equivalent, simplified form of this result (in which $ A_j = A_{j0}$, $B_j = S-A_{jq}$, and $n$ replaces $n-q$) is

\begin{theorem}\mylabel{gst} 
Let $n > 0$.  Suppose $ \left( A_j , B_j \right) $ are $m$ pairs of sets such that 
$ A_j \cap B_j = \eset $ for all $j$, $ A_j \cap B_k \neq \eset $ for all $j \neq k$, and all 
$ \left| A_j \right| + \left| B_j \right| \leq n $. Then $ m \leq \binom{n}{ \lfloor n/2 \rfloor } $ and this bound can be attained for every $n$. \end{theorem} 

Sperner's inequality follows as the special case in which $ A_1, \dots , A_m \subseteq S $ form an antichain and $ B_j = S - A_j $. 

Theorems \ref{erdos}, \ref{meshalkin}, and \ref{gst} are \emph{incomparable} generalizations of Sperner's Theorem. 
We wish to combine (and hence further generalize) these generalizations. 
To state our main result, we define a {\bf weak partial composition of $S$ into $p$ parts} as an ordered $p$-tuple 
$ \left( A_1, \dots , A_p \right) $ such that $ A_1, \dots , A_p $ are pairwise disjoint sets, possibly void (hence the word ``weak''), and $ A_1 \cup \dots \cup A_p \subseteq S $. 
If we do not specify the superset $S$ then we simply talk about a {\bf weak set composition into $p$ parts} (this could be a weak composition of any set). 
Our generalization of Sperner's inequality is:
\begin{theorem}\mylabel{grandmother} 
Fix integers $p \geq 2$ and $r \geq 1$.  Suppose $ \left( A_{j1} , \dots , A_{jp} \right) $ for $j=1,\ldots,m$ are different weak set compositions 
into $p$ parts with the condition that, for all $ k \in [p] $ and all $ I \subseteq [m] $ with $ |I| = r+1 $, there exist distinct 
$ i, j \in I$ such that either $ A_{ik} = A_{jk} $ or 
\begin{equation}\mylabel{condition} 
A_{ik} \cap \bigcup_{ l \neq k } A_{jl} \neq \eset \neq A_{jk} \cap \bigcup_{ l \neq k } A_{il} \ , 
\end{equation} 
and let $ n := \max_{ 1 \leq j \leq m } \left( | A_{j1} | + \dots + | A_{jp} | \right) $. 
Then $m$ is bounded by the sum of the $ r^p $ largest $p$-multinomial coefficients for 
integers less than or equal to $n$. \end{theorem}

If $ r^p $ is larger than $ \binom{r+p}{p} $, the number of $p$-multinomial coefficients, then we regard the sequence of coefficients as extended by $0$'s.

We heartily agree with those readers who find the statement of this theorem somewhat unreadable. 
We would first like to show that it does generalize Theorems \ref{erdos}, \ref{meshalkin}, and \ref{gst} simultaneously. 
The last follows easily as the case $ r=1 $, $ p=2 $. Theorem \ref{meshalkin} can be deduced by choosing $r=1$ and 
restricting the weak compositions to be compositions of a fixed set $S$ with $n$ elements. Finally, Theorem \ref{erdos} 
follows by choosing $p=2$ and the weak compositions to be compositions of a fixed $n$-set into 2 parts. 

What we find more interesting, however, is that specializations of Theorem \ref{grandmother} yield simply stated corollaries that combine 
two at a time of Theorems \ref{erdos}, \ref{meshalkin}, and \ref{gst}. Section \ref{concl} collects these corollaries. 

The condition of the theorem implies that each set $\cA_k = \{ A_{jk} : j \in [m] \}$ (ignoring repetition) is $r$-chain-free. 
We suspect that the converse is not true in general.  (It is true if all the weak set compositions are weak compositions of the same set 
of order $n$, as in Corollary \ref{e-m}.) 

All the theorems we have stated have each a slightly stronger companion, an {\em LYM inequality}. 
In Section \ref{lym}, we state these inequalities and show how Theorems \ref{sperner}--\ref{grandmother} 
can be deduced from them. The proofs of Theorem \ref{grandmother} and the corresponding LYM inequality are in Section \ref{proof}. 
After the corollaries of Section \ref{concl}, in Section \ref{max} we show that some, at least, of our upper bounds cannot be attained. 


\section{LYM inequalities}\mylabel{lym} 
In attempting to find a new proof of Theorem \ref{sperner}, Lubell, Yamamoto, and Meshalkin independently came up with the following refinement: 

\begin{theorem}[Lubell \cite{lubell}, Yamamoto \cite{yamamoto}, Meshalkin \cite{meshalkin}]\mylabel{lymineq} Suppose $ A_1, \dots , A_m \subseteq S $ such that $ A_k \not\subseteq A_j $ for $ k \neq j $. Then 
\[ 
\sum_{ k=1 }^{ m } \frac{1}{ \binom{ n }{ \left| A_k \right| } } \leq 1 \ . 
\] 
\end{theorem}

Sperner's inequality follows immediately by noting that 
$ \max_k \binom{n}{k} = \binom{n}{ \lfloor n/2 \rfloor } $ . 

An LYM inequality corresponding to Theorem \ref{erdos} appeared to our knowledge first in \cite{harper}: 

\begin{theorem}[Rota--Harper]\mylabel{erdoslym} Suppose $ \left\{ A_1, \dots , A_m \right\} \subseteq \cP(S) $ contains no chains with $r+1$ elements. Then 
\[ 
\sum_{ k=1 }^{ m } \frac{1}{ \binom{ n }{ \left| A_k \right| } } \leq r \ . 
\] 
\end{theorem}

Deducing Erd\H{o}s's Theorem \ref{erdos} from this inequality is not as straightforward as the connection between Theorems \ref{lymineq} and \ref{sperner}. 
It can be done through Lemma \ref{cute}, which we also need in order to deduce Theorem \ref{grandmother}. 

The LYM companion of Theorem \ref{meshalkin} first appeared in \cite{hochberg}; again, Meshalkin's Theorem \ref{meshalkin} follows immediately. 

\begin{theorem}[Hochberg--Hirsch]\mylabel{meshalkinlym} Suppose $ \left( A_{j1} , \dots , A_{jp} \right) $ for $j=1,\ldots,m$ are different weak compositions of $S$ into $p$ parts 
such that for each $ k \in [p] $ the set $ \left\{ A_{jk} : 1 \leq j \leq m \right\} $ (ignoring repetitions) forms 
an antichain. Then 
\[ 
\sum_{ j=1 }^{ m } \frac{1}{ \binom{ n }{ \left| A_{j1} \right| , \dots , \left| A_{jp} \right| } } \leq 1 \ . 
\] 
\end{theorem} 

The LYM inequality corresponding to Theorem \ref{gst} is due to Bollob\'as. 

\begin{theorem}[Bollob\'as \cite{bollobas}]\mylabel{bollobas} Suppose $ \left( A_j , B_j \right) $ are $m$ pairs of sets such that 
$ A_j \cap B_j = \eset $ for all $j$ and $ A_j \cap B_k \neq \eset $ for all $j \neq k$. Then 
\[ 
\sum_{ j=1 }^{ m } \frac{1}{ \binom{ \left| A_{j} \right| + \left| B_{j} \right| }{ \left| A_{j} \right| } } \leq 1 \ . 
\] 
\end{theorem} 

Once more, the corresponding upper bound, the Griggs--Stahl--Trotter Theorem \ref{gst}, is an immediate consequence. 

Naturally, there is an LYM inequality accompanying our main Theorem \ref{grandmother}. Like its siblings, it constitutes a refinement.

\begin{theorem}\mylabel{grandmotherlym} 
Let $p \geq 2$ and $r \geq 1$.  Suppose $ \left( A_{j1} , \dots , A_{jp} \right) $ for $j=1,\ldots,m$ are different weak compositions (of any sets) 
into $p$ parts satisfying the same condition as in Theorem \ref{grandmother}.
Then 
\[ 
\sum_{j=1}^{m} \frac{1}{ \binom{ | A_{j1} | + \dots + | A_{jp} | }{ | A_{j1} | , \dots , | A_{jp} | } }  \leq r^p \ . 
\] 
\end{theorem} 

\comment{Added August, 2001.}
\begin{example}\mylabel{notr} 
The complicated hypothesis of Theorem \ref{grandmotherlym} cannot be replaced by the assumption that each $\cA_k$ is $r$-chain-free, because then there is no LYM bound independent of $n$.  
Let $n \gg p \geq 2$, $S=[n]$, and $\cA = \{ (A,\{n\},\{n-1\},\ldots,\{n-p+2\}) : A \in \cA_1 \}$ where $\cA_1$ is a largest $r$-chain-free family in $[n-p+1]$, specifically, 
$$
\cA_1 = \bigcup_{j\in I} \cP_j\big([n-p+1]\big) 
$$
\nopagebreak
where
$$
I = 
\left\{ \left\lceil\frac{n-p+1-r}{2}\right\rceil, \left\lceil\frac{n-p+1-r}{2}\right\rceil +1, \left\lceil\frac{n-p+1-r}{2}\right\rceil +r-1 \right\}.
$$
The LYM sum is
\begin{align*}
\sum_{A\in\cA_1} \frac{1}{\binom{|A|+p-1}{|A|,1,\ldots,1}} 
&= \sum_{A\in\cA_1} \frac{|A|!}{(|A|+p-1)!} 	\\
&= \sum_{j\in I} \binom{n-p+1}{j} \frac{j!}{(j+p-1)!} \\
&= \sum_{j\in I} \frac{(n-p+1)\cdots(n-p-j+2)}{(p-1+j)!} \\
&\to \infty \quad\text{as } n\to \infty.
\end{align*}
There is no possible upper bound in terms of $n$.
\end{example}


\section{Proof of the main theorems}\mylabel{proof} 
{\it Proof of Theorem} \ref{grandmotherlym}. Let $S$ be a finite set containing all $ A_{jk} $ for $ j = 1, \dots, m $ and $ \ k = 1, \dots, p $, and let $ n = |S| $. We count maximal chains in $ \cP(S) $. 

Let us say a maximal chain {\bf separates} the weak composition $ \left( A_{1} , \dots , A_{p} \right) $ if there exist elements 
$ \eset = X_0 \subseteq X_{l_1} \subseteq \dots \subseteq X_{l_p} = S $ of the maximal chain such that 
$ A_k \subseteq X_{l_k} - X_{l_{k-1}} $ for each $k$. 
There are 
\begin{equation} \mylabel{star}
\binom{ n }{ | A_{1} | + \cdots + | A_{p} | } | A_{1} | ! \cdots | A_{p} | ! \left( n - | A_{1} | - \dots - | A_{p} | \right) ! 
\end{equation}
maximal chains separating $ \left( A_{1} , \dots , A_{p} \right) $.  (To prove this, replace maximal chains $\eset \subset \{x_1\} \subset \{x_1,x_2\} \subset \cdots \subset S$ by permutations $(x_1,x_2,\ldots,x_n)$ of $S$.  Choose $| A_{1} | + \cdots + | A_{p} |$ places for $A_1 \cup \cdots \cup A_p$; then arrange $A_1$ in any order in the first $|A_1|$ of these places, $A_2$ in the next $|A_2|$, etc.  Finally, arrange $S - ( A_1 \cup \cdots \cup A_p)$ in the remaining places.  This constructs all maximal chains that separate $ \left( A_{1} , \dots , A_{p} \right) $.)

We claim that every maximal chain separates at most $ r^p $ weak partial compositions of $|S|$. 
To prove this, assume that there is a maximal chain that separates $ N $ weak partial compositions $ \left( A_{j1} , \dots , A_{jp} \right) $. 
Consider all first components $ A_{j1} $ and suppose $ r+1 $ of them are different, say $ A_{11}, A_{21}, \ldots, A_{r+1,1} $.  By the hypotheses of the theorem, there are $i, i' \in [r+1]$ such that $A_{i1}$ meets some $A_{i'l'}$ where $l' > 1$ and $A_{i'1}$ meets some $ A_{il} $ where $ l>1 $.  By separation, there are $q_1^{}$ and $q'_1$ such that $ A_{i1} \subseteq X_{q_1} - X_0 $ and $ A_{i'1} \subseteq X_{q'_1} - X_0 $, and there are $ q_{l-1}^{}, q_l^{}, q'_{l'-1}, q'_{l'} $ such that $ q_1^{} \leq q_{l-1}^{} \leq q_l^{} $, $ q'_1 \leq q'_{l'-1} \leq q'_{l'} $, and 
\[
A_{il} \subseteq X_{q_l^{}} - X_{q_{l-1}^{}} \qquad \text{ and } \qquad 
A_{i'l'} \subseteq X_{q'_{l'}} - X_{q'_{l'-1}} \ .
\]
Since $ A_{i1} $ meets $ A_{i'l'} $, there is an element $ a_{i1} \in X_{q'_{l'}} - X_{q'_{l'-1}} $; it follows that $ q'_{l'-1} < q_1^{} $.  Similarly, $ q_{l-1}^{} < q'_1 $.  But this is a contradiction.  It follows that, amongst the $N$ sets $A_{j1}$, there are at most $r$ different sets.  Hence (by the pigeonhole principle) there are $ \lceil N/r \rceil $ among the $N$ weak partial compositions that have the same first set $A_{j1}$.

Looking now at these $ \lceil N/r \rceil $ weak partial compositions, we can repeat the argument to conclude that there are $ \big\lceil \lceil N/r \rceil / r \big\rceil \geq \left\lceil N/r^2 \right\rceil $ weak partial compositions for which both the $ A_{j1} $'s and the $ A_{j2} $'s are identical. 
Repeating this process $p-1$ times yields $ \left\lceil N/r^{p-1} \right\rceil $ weak partial compositions into $p$ parts whose first $p-1$ parts 
are identical. But now the hypotheses imply that the last parts of all these weak partial compositions are at most $r$ different sets; 
in other words, there are at most $r$ distinct weak partial compositions. Hence $ \left\lceil N/r^{p-1} \right\rceil \leq r $, whence $ N \leq r^p $. 
(If we know that all the compositions are weak---but not partial---compositions of $S$, then the last parts of all these $ \left\lceil N/r^{p-1} \right\rceil $ 
weak compositions are identical. Thus $ N \leq r^{p-1} $.) 

Since at most $r^p$ weak partial compositions of $S$ are separated by each of the $n!$ maximal chains, from (\ref{star}) we deduce that
\begin{eqnarray*}
r^p n! \geq \sum_{j=1}^{m} \binom{ n }{ | A_{j1} | + \dots + | A_{jp} | } | A_{j1} | ! \cdots | A_{jp} | ! \left( n - | A_{j1} | - \dots - | A_{jp} | \right) ! \\ 
= \sum_{j=1}^{m} \frac{n!}{ \binom{ |A_{j1}| + \dots + |A_{jp}| }{ | A_{j1} | , \dots , | A_{jp} | } }  \ . 
\end{eqnarray*}
The theorem follows.
\hfill {} $\qedsymbol$

To deduce Theorem \ref{grandmother} from Theorem \ref{grandmotherlym}, we use the following lemma, which originally appeared in somewhat different and incomplete form in \cite{harper}, used there to prove Erd\H{o}s's Theorem \ref{erdos} 
by means of Theorem \ref{erdoslym}, and appeared in complete form in \cite[Lemma 3.1.3]{klain}. We give a very short proof, which seems to be new. 

\begin{lemma}[Harper--Klain--Rota] \mylabel{cute} 
Suppose $ M_1, \dots, M_N \in \bbR $ satisfy $ M_1 \geq M_2 \geq \dots \geq M_N \geq 0$,
and let $R$ be an integer with $ 1 \leq R \leq N $. If $ q_1, \dots, q_N \in [0,1] $ have sum
\[ 
q_1 + \dots + q_N \leq R \ ,
\]
then 
\[
q_1 M_1 + \dots + q_N M_N \leq M_1 + \dots + M_R \ . 
\] 
\end{lemma}

{\it Proof.} By assumption,
\[ 
\sum_{ k=R+1 }^{ N } q_k \ \leq \ \sum_{ k=1 }^{ R } ( 1 - q_k ) \ . 
\]
Hence, by the condition on the $ M_k $,
\[ 
\sum_{ k=R+1 }^{ N } q_k M_k \ \leq \ M_R \sum_{ k=R+1 }^{ N } q_k \ \leq \ M_R \sum_{ k=1 }^{ R } ( 1 - q_k ) \ \leq \ \sum_{ k=1 }^{ R } ( 1 - q_k ) M_k \ , 
\] 
which is equivalent to the conclusion.
\hfill {} $\qedsymbol$

{\it Proof of Theorem} \ref{grandmother}.  
Let $S$ be any finite set that contains all $A_{jk}$.  
Write down the LYM inequality from Theorem \ref{grandmotherlym}.  


From the $m$ weak partial compositions $ \left( A_{j1} , \dots , A_{jp} \right) $ of $S$, 
collect those whose shape is $ ( a_1, \dots, a_p ) $ into the set $ C( a_1 , \dots , a_p ) $.  
Label the $p$-multinomial coefficients for integers $n' \leq n$ as $M_1', M_2', \ldots$ so that $ M_1' \geq M_2' \geq \cdots $. 
If $ M_k' $ is $ \binom{ n' }{ a_1, \dots, a_p } $, 
let $ q_k' := \left| C( a_1 , \dots , a_p ) \right| / M_k' $. By Theorem \ref{grandmotherlym}, the $ q_k' $'s and $ M_k'  $'s satisfy all the conditions of 
Lemma \ref{cute} with $N$ replaced by the number of $p$-tuples $ ( a_1 , \dots , a_p ) $ whose sum is at most $n$, 
that is $ \binom{n+p}{p} $, and $R$ replaced by $ \min( N, r^p ) $.  Hence 
\[ 
\sum_{ a_1 + \dots + a_p \leq n } \left| C( a_1 , \dots , a_p ) \right| \leq M_1' + \dots + M_{R}' \ . 
\] 
The conclusion of the theorem now follows, since 
\[ 
m = \sum_{ a_1 + \dots + a_p \leq n } \left| C( a_1 , \dots , a_p ) \right| \ . 
\] 
\hfill {} $\qedsymbol$


\section{Consequences}\mylabel{concl} 
As promised in Section \ref{intro}, we now state special cases of Theorems \ref{grandmother}/\ref{grandmotherlym} 
that unify pairs of Theorems \ref{erdos}, \ref{meshalkin}, and \ref{gst} as well as their LYM companions. 

The first special case unifies Theorems \ref{erdos}/\ref{erdoslym} and \ref{meshalkin}/\ref{meshalkinlym}. 
(It is a corollary of the proof of the main theorems, not of the theorems themselves.  See \cite{bz} for a short, direct proof.) 

\begin{corollary} \mylabel{e-m}
Suppose $ \left( A_{j1} , \dots , A_{jp} \right) $ are $m$ different weak compositions of $S$ 
into $p$ parts such that for each $ k \in [p-1] $, 
the set $ \left\{ A_{jk} : 1 \leq j \leq m \right\} $ is $r$-chain-free. Then 
\[ 
\sum_{j=1}^{m} \frac{1}{ \binom{ n }{ | A_{j1} | , \dots , | A_{jp} | } }  \leq r^{p-1} \ . 
\] 
Consequently, $m$ is bounded by the sum of the $ r^{p-1} $ largest $p$-multinomial coefficients for $n$. 
\end{corollary} 

{\it Proof.} We note that, for a family of $m$ weak compositions of $S$, the condition of Theorem \ref{grandmotherlym} 
for a particular $k \in [p-1]$ is equivalent to $ \{ A_{jk} \}_j $ being $r$-chain-free. Thus by the hypothesis of the 
corollary, the hypothesis of the theorem is met for $k=1,\dots,p-1$. Then the proof of Theorem \ref{grandmotherlym} 
goes through perfectly with the only difference, explained in the proof, that (even without a condition on $k=p$) 
we obtain $ N \leq r^{p-1} $. In the proof of Theorem \ref{grandmother}, under our hypotheses the sets $C(a_1,\dots,a_p)$ 
with $a_1+\dots+a_p < n$ are empty. Therefore we take only the $p$-multinomial coefficients for $n$, labelled 
$M_1 \geq M_2 \geq \cdots$. In applying Lemma \ref{cute} we take $R = \min (N,r^{p-1})$ and summations over 
$a_1+\dots+a_p = n$. With these alterations the proof fits Corollary \ref{e-m}. 
\hfill {} $\qedsymbol$

\comment{Corollary \ref{rfamily} and discussion added Aug.\ 2001:}
A good way to think of Corollary \ref{e-m} is as a theorem about partial weak compositions, obtained by dropping the last part from each of the weak compositions in the corollary.

\begin{corollary} \mylabel{rfamily}
Fix $p \geq 2$ and $r \geq 1$.  Suppose $ \left( A_{j1} , \dots , A_{jp} \right) $ are $m$ different weak partial compositions 
of an $n$-set $S$ into $p$ parts such that for each $ k \in [p] $, 
the set $ \left\{ A_{jk} : 1 \leq j \leq m \right\} $ is $r$-chain-free. Then 
$m$ is bounded by the sum of the $ r^{p} $ largest $(p+1)$-multinomial coefficients for $n$. 
\hfill $\qedsymbol$
\end{corollary} 

A difference between this and Theorem \ref{grandmother} is that Corollary \ref{rfamily} has a weaker and simpler hypothesis but a much weaker bound.  
But the biggest difference is the omission of an accompanying LYM inequality.  Corollary \ref{e-m} obviously implies one, but it is weaker than that in Theorem \ref{grandmotherlym} because, since the top number in the latter can be less than $n$, the denominators are much smaller.  We do not present in Corollary \ref{rfamily} an LYM inequality of the kind in Theorem \ref{grandmotherlym} for the very good reason that none is possible; that is the meaning of Example \ref{notr}.

\comment{
We would get a stronger result by dropping the last part from the weak compositions, 
forming $ A_i' = ( A_{i,1} , \dots , A_{i,p-1} ) $ and applying the main theorems to $A'_1,\ldots,A'_m$. However, we 
cannot see how to deduce the overlap hypotheses of these theorems from the $r$-family hypotheses of Corollary \ref{e-m}. 
}

The second specialization constitutes a weak common refinement of Theorems 
\ref{erdos}/\ref{erdoslym} and \ref{gst}/\ref{bollobas}.  We call it weak because its specialization to the case $B_j = S - A_j $, which is the situation of Theorems \ref{erdos}/\ref{erdoslym}, is weaker than those theorems. 

\begin{corollary} \mylabel{e-g}
Let $r$ be a positive integer.  Suppose $ \left( A_j , B_j \right) $ are $m$ pairs of sets such that 
$ A_j \cap B_j = \eset $ and, 
for all $ I \subseteq [m] $ with $ |I| = r+1 $, there exist distinct $ i, j \in I $ 
for which $ A_j \cap B_k \neq \eset \neq A_k \cap B_j $.
Let $ n = \max_j \left( \left| A_j \right| + \left| B_j \right| \right)$.  Then 
\[ 
\sum_{ j=1 }^{ m } \frac{1}{ \binom{ \left| A_{j} \right| + \left| B_{j} \right| }{ \left| A_{j} \right| } } \leq r \ . 
\] 
Consequently, $m$ is bounded by the sum of the $r$ largest binomial coefficients $ \binom{n'}{k} $ for $ 0 \leq k \leq n' \leq n $. 
This bound can be attained for all $n$ and $r$. 
\end{corollary}  

{\it Proof.}  Set $p=2$ in Theorems \ref{grandmother}/\ref{grandmotherlym}.  To attain the bound, let $A_j$ range over all $k$-subsets of $[n]$ and let $B_j=[n]-A_j$.  
\hfill $\qedsymbol$

The last special case of Theorems \ref{grandmother}/\ref{grandmotherlym} we would like to mention is that in which $ r=1 $; it unifies 
Theorems \ref{meshalkin}/\ref{meshalkinlym} and \ref{gst}/\ref{bollobas}.

\begin{corollary} \mylabel{m-g}
Suppose $ \left( A_{j1} , \dots , A_{jp} \right) $ are $m$ different weak set compositions 
into $p$ parts with the condition that, for all $ k \in [p] $ and all distinct 
$ i, j \in [m] $, either $ A_{ik} = A_{jk} $ or 
\[ 
A_{ik} \cap \bigcup_{ l \neq k } A_{jl} \neq \eset \neq A_{jk} \cap \bigcup_{ l \neq k } A_{il} \ . 
\] 
and let $ n \geq \max_j \big( | A_{j1} | + \dots + | A_{jp} | \big) $.  Then 
\[ 
\sum_{j=1}^{m} \frac{1}{ \binom{ | A_{j1} | + \dots + | A_{jp} | }{ | A_{j1} | , \dots , | A_{jp} | } }  \leq 1 \ . 
\] 
Consequently, $m$ is bounded by the largest $p$-multinomial coefficient for $n$. 
The bound can be attained for every $n$ and $p$.
\end{corollary} 

{\it Proof.}  Everything follows from Theorems \ref{grandmother}/\ref{grandmotherlym} except the attainability of the upper bound, which is a consequence of Theorem \ref{meshalkin}.
\hfill $\qedsymbol$ 


\section{The maximum number of compositions} \mylabel{max} 
Although the bounds in all the previously known Sperner generalizations of Section \ref{intro} can be attained, 
for the most part that seems not to be the case in Theorem \ref{grandmother}. 
The key difficulty appears in the combination of $r$-families with compositions as in Corollary \ref{e-m}. 
(We think it makes no difference if we allow partial compositions but we have not proved it.)  We begin with a refinement of Lemma \ref{cute}. 
A weak set composition has {\bf shape} $ (a_1,\dots,a_p) $ if $ | A_k | = a_k $ for all $k$. 

%
%
\begin{lemma} \mylabel{sharp} 
Given values of $n$, $r$, and $p$ such that $ r^{p-1} \leq \binom{ n+p-1 }{ p-1 } $, the bound in Corollary \ref{e-m} 
can be attained only by taking all weak compositions of shape $ ( a_1,\ldots,a_p ) $ that give $p$-multinomial coefficient larger than the 
$ (r^{p-1}+1) $-st largest such coefficient $M_{r^{p-1}+1}$, and none whose shape gives a smaller coefficient than 
the $ (r^{p-1}) $-st largest such coefficient $M_{r^{p-1}}$. 
\end{lemma} 
{\it Proof.}  
First we need to characterized sharpness in Lemma \ref{cute}.  Our lemma is a slight improvement on \cite[Lemma 3.1.3]{klain}.

\begin{lemma} \mylabel{cutesharp}
In Lemma \ref{cute}, suppose that $M_R > 0$. 
Then there is equality in the conclusion if and only if
\[
q_k^{} = 1 \text{ if } M_k > M_R \qquad\text{ and }\qquad q_k^{} = 0 \text{ if } M_k < M_R 
\]
and also, letting $M_{R'+1}$ and $M_{R''}$ be the first and last $M_k$'s equal to $M_R$, 
\[
q_{R'+1}^{} + \dots + q_{R''}^{} = R - R' \ .
\qquad\qquad \qedsymbol
\]
\end{lemma}

In Lemma \ref{sharp}, all $M_k > 0$ for $ k \leq \binom{ n+p-1 }{ p-1 } $. 
(We assume $N$ is no larger than $ \binom{ n+p-1 }{ p-1 } $.  The contrary case is easily derived from that one.)
It is clear that, when applying Lemma \ref{cute}, we have to have in our set of weak compositions all those of the shapes 
$ ( a_1,\ldots,a_p ) $ for which $ \binom{ n }{ a_1, \dots, a_p } > M_{r^{p-1}} $ and none for which $ \binom{ n }{ a_1, \dots, a_p } < M_{r^{p-1}} $.  
The rest of the $m$ weak compositions can have any shapes for which $ \binom{ n }{ a_1, \dots, a_p } = M_{r^{p-1}} $. 
If $ M_{r^{p-1}} > M_{r^{p-1}+1} $ this means we must have all weak compositions with shapes for which $ \binom{ n }{ a_1, \dots, a_p } > M_{r^{p-1}+1} $. 
\hfill $\qedsymbol$

To explain why the bound cannot usually be attained, we need to define the ``first appearance'' of a size $ a_i $ in the descending 
order of $p$-multinomial coefficients for $n$. 

Fix $ p \geq 3 $ and $n$ and let $ n = \nu p + \rho $ where $ 0 \leq \rho < p $. In $ \binom{ n }{ a_1, \dots, a_p } $, the $ a_i $ are the 
{\bf sizes}. The multiset of sizes is the {\bf form} of the coefficient. Arrange the multinomial coefficients in decreasing order: 
$ M_1 \geq M_2 \geq M_3 \geq \cdots $. (There are many such orderings; choose one arbitrarily, fix it, and call it 
{\bf the descending order} of coefficients.) Thus, for example, 
\[ 
M_1 = \binom{ n }{ \nu, \dots, \nu } > M_2 = \binom{ n }{ \nu + 1, \nu, \dots, \nu, \nu - 1 } = M_3 = \dots = M_{p(p-1)+1} \qquad \mbox{ if } p | n 
\] 
since $ M_3, \dots , M_{p(p-1)+1} $ have the same form as $ M_2 $, and 
\[ 
M_1 = \binom{ n }{ \nu + 1, \dots, \nu } = \dots = M_{ \binom p \rho } > M_{ \binom p \rho + 1 } \qquad \mbox{ if } p \nm n , 
\] 
where the form of $ M_1 $ has $ \rho $ sizes equal to $ \nu + 1 $, so $ M_1 , \dots , M_{ \binom p \rho } $ all have the same form. 

As we scan the descending order of multinomial coefficients, each possible size $ \kappa, 0 \leq \kappa \leq n $, appears first 
in a certain $ M_i $. We call $ M_i $ the {\bf first appearance} of $ \kappa $ and label it $ L_\kappa $. 
For example, if $ p | n $, $ L_\nu = M_1 > L_{ \nu + 1 } = L_{ \nu - 1 } = M_2 $, while if $ p \nm n $ then $ L_\nu = L_{ \nu + 1 } = M_1 $. 
It is clear that $ L_\nu > L_{ \nu - 1 } > \dots $ and $ L_{ \nu + 1 } > L_{ \nu + 2 } > \dots $, but the way in which the lower $ L_\kappa $'s, 
where $ \kappa \leq \nu $, interleave the upper ones is not obvious. 
We write $ L_k^\s $ for the $k$-th $ L_\kappa $ in the descending order of multinomial coefficients. Thus $ L_1^\s = L_\nu $; 
$ L_2^\s = L_{\nu+1} $ and $ L_3^\s = L_{\nu-1} $ (or vice versa) if $ p | n $, and $ L_2^\s = L_{\nu+1} $ if $ p \nm n $ 
while $ L_3^\s = L_{\nu+2} $ or $ L_{\nu-1} $. 
\begin{theorem}\mylabel{attain} 
Given $ r \geq 2 $, $ p \geq 3 $, and $ n \geq p $, the bound in Corollary \ref{e-m} cannot be attained if $ L_r^\s > M_{ r^{ p-1 } + 1 } $. 
\end{theorem} 
The proof depends on the following lemma. 
\begin{lemma}\mylabel{attainlemma} 
Let $ r \geq 2 $ and $ p \geq 3 $, and let $ \kappa_1 , \dots , \kappa_r $ be the first $r$ sizes that appear in the descending order of 
$p$-multinomial coefficients for $n$. The number of all coefficients with sizes drawn from $ \kappa_1 , \dots , \kappa_r $ is less than 
$ r^{p-1} $ and their sum is less than $ M_1 + \dots + M_{r^{p-1}} $. 
\end{lemma} 
{\it Proof.} Clearly, $ \kappa_1 , \dots , \kappa_r $ form a consecutive set that includes $ \nu $. Let $ \kappa $ be the smallest and 
$ \kappa' $ the largest. One can verify that, in $ \binom n { \kappa , \dots , \kappa , x } $ and $ \binom n { \kappa' , \dots , \kappa' , y } $, 
it is impossible for both $x$ and $y$ to lie in the interval $ [ \kappa , \kappa' ] $ as long as $ (r-1)(p-2) > 0 $. 
\hfill $\qedsymbol$ 

{\it Proof of Theorem} \ref{attain}. Suppose the upper bound of Corollary \ref{e-m} is attained by a certain set of weak compositions of $S$, 
an $n$-element set. For each of the first $r$ sizes $ \kappa_1 , \dots , \kappa_r $ that appear in the descending order of $p$-multinomial 
coefficients, $ L_{\kappa_i} $ has sizes drawn from $ \kappa_1 , \dots , \kappa_r $ and at least one size $ \kappa_i $. 
Taking all coefficients $ M_k $ that have the same forms as the $L_{\kappa_i} $, $ \kappa_i $ will appear in each position $j$ in some $M_k$. 
By hypothesis and Lemma \ref{sharp}, among our set of weak compositions, every $ \kappa_i $-subset of $S$ appears in every position in the 
weak compositions. If any subset of $S$ of a different size from $ \kappa_1 , \dots , \kappa_r $ appeared in any position, there would be a 
chain of length $r$ in that position. Therefore we can only have weak compositions whose sizes are among the first $r$ sizes. 
By Lemma \ref{attainlemma}, there are not enough of these to attain the upper bound. 
\hfill $\qedsymbol$ 

Theorem \ref{attain} can be hard to apply because we do not know $ M_{ r^{p-1} + 1 } $. On the other hand, we do know $ L_\kappa $ since it 
equals $ \binom n { \kappa, a_2 , \dots, a_p } $ where $ a_2, \dots, a_p $ are as nearly equal as possible. A more practical criterion 
for nonattainment of the upper bound is therefore 
\begin{corollary}\mylabel{attaincor} 
Given $ r \geq 2 $, $ p \geq 3 $, and $ n \geq p $, the bound in Corollary \ref{e-m} cannot be attained if $ L_r^\s > L_{r+1}^\s $. 
\end{corollary} 
{\it Proof.} It follows from Lemma \ref{attainlemma} that $ L_{r+1}^\s $ is one of the first $ r^{p-1} $ coefficients. Thus 
$ L_r^\s > L_{r+1}^\s \geq M_{ r^{p-1} +1 } $ and Theorem \ref{attain} applies. 
\hfill $\qedsymbol$ 

It seems clear that $ L_r^\s $ will almost always be larger than $ L_{r+1}^\s $ (if $ r \geq 3 $ or $ p \nm n $) so our bound 
will not be attained. However, cases of equality do exist. For instance, take $p=3$, $ r=3 $, and $n=10$; then 
$ L_5^\s = L_1 = \binom{10}{5,4,1} = 1260 $ and $ L_6^\s = L_6 = \binom{10}{6,2,2} = 1260 $. Thus if $r=5$, Corollary \ref{attaincor} 
does not apply here. (We think the bound is still not attained but we cannot prove it.) 
We can isolate the instances of equality for each $r$, but as $r$ grows larger the calculations quickly become extensive. 
Thus we state the results only for small values of $r$. 
\begin{proposition} 
The bound in Corollary \ref{e-m} cannot be attained if $ 2 \leq r \leq 5 $ and $ p \geq 3 $ and $ n \geq r-1 $, 
except possibly when $ r=2 $, $ p | n $, and $ p=3,4,5$, or when $ r=4 $, $ p \geq 4 $, and $ n=2p-1 $, or when $r=5$, $p=3$, and $n=10$. 
\end{proposition} 
{\it Proof sketch.} Suppose $ p \nm n $. We have verified (by long but routine calculations which we omit) that 
$ L_1^\s = L_2^\s > L_3^\s > L_4^\s > L_5^\s > L_6^\s $ except that $ L_4^\s = L_5^\s $ if $ \rho = p-1 $ and $ p \geq 4 $ and $ \nu = 1 $ 
and $ L_5^\s = L_6^\s $ when $p=\nu=3$ and $\rho=1$. 

If $ p | n $ then $ L_1^\s > L_2^\s = L_3^\s > L_4^\s > L_5^\s > L_6^\s $. This implies the proposition for $r=3$, $4$, or $5$. 
We approach $ r=2 $ differently. The largest coefficients are 
\[ 
M_1 = \binom n { \nu, \dots, \nu } > M_2 = \binom n { \nu + 1, \nu, \dots, \nu, \nu - 1 } = \dots = M_{ p(p-1) + 1 } > M_{ p(p-1) + 2 } \ . 
\] 
If $ p(p-1) + 1 \leq r^{p-1} $, the bound is unattainable by Theorem \ref{attain}. That is the case when $ p \geq 6 $. 




\small
\nocite{*}
\bibliography{thesis}
\bibliographystyle{alpha}

\end{document}